\documentclass[11pt]{article}

\usepackage{amsmath}
\usepackage{amsthm} 
\usepackage{amssymb} 
\usepackage{hyperref}
\usepackage{color} 

\usepackage{enumitem}

\newtheorem{thm}{Theorem}
\newtheorem{inspr}[thm]{}

\newenvironment{definitie}{\begin{itemize}\item[ ]\hspace{-26pt}\bf Definition \rm }{\end{itemize}}
\newenvironment{notatie}{\begin{itemize}\item[ ]\hspace{-26pt}\bf Notation \rm }{\end{itemize}}
\newenvironment{voorbeeld}{\begin{itemize}\item[ ]\hspace{-26pt}\bf Example \rm }{\end{itemize}}
\newenvironment{stelling}{\begin{itemize}\item[ ]\hspace{-26pt}\bf Theorem \rm }{\end{itemize}}
\newenvironment{propositie}{\begin{itemize}\item[ ]\hspace{-26pt}\bf Proposition \rm }{\end{itemize}}
\newenvironment{lemma}{\begin{itemize}\item[ ]\hspace{-26pt}\bf Lemma \rm }{\end{itemize}}
\newenvironment{opmerking}{\begin{itemize}\item[ ]\hspace{-26pt}\bf Remark \rm }{\end{itemize}}
\newenvironment{voorwaarde}{\begin{itemize}\item[ ]\hspace{-26pt}\bf Condition \rm }{\end{itemize}}
\newenvironment{probleem}{\begin{itemize}\item[ ]\hspace{-26pt}\bf Problem \rm }{\end{itemize}}

\newcommand{\defin}{\begin{inspr}\begin{definitie}}  
\newcommand{\edefin}{\end{definitie}\end{inspr}}

\newcommand{\notat}{\begin{inspr}\begin{notatie}}  
\newcommand{\enotat}{\end{notatie}\end{inspr}}

\newcommand{\voorb}{\begin{inspr}\begin{voorbeeld}}  
\newcommand{\evoorb}{\end{voorbeeld}\end{inspr}}

\newcommand{\stel}{\begin{inspr}\begin{stelling}}
\newcommand{\estel}{\end{stelling}\end{inspr}}

\newcommand{\prop}{\begin{inspr}\begin{propositie}}
\newcommand{\eprop}{\end{propositie}\end{inspr}}

\newcommand{\lem}{\begin{inspr}\begin{lemma}}
\newcommand{\elem}{\end{lemma}\end{inspr}}

\newcommand{\opm}{\begin{inspr}\begin{opmerking}}
\newcommand{\eopm}{\end{opmerking}\end{inspr}}

\newcommand{\voorw}{\begin{inspr}\begin{voorwaarde}}
\newcommand{\evoorw}{\end{voorwaarde}\end{inspr}}

\newcommand{\probl}{\begin{inspr}\begin{probleem}}
\newcommand{\eprobl}{\end{probleem}\end{inspr}}

\newcommand{\bew}{\vspace{-0.3cm}\begin{itemize}\item[ ] \bf Proof\rm: }
\newcommand{\bewo}{\vspace{-0.3cm}\begin{itemize}\item[ ] \bf Proof\rm (outline): }
\newcommand{\ebew}{\hfill $\qed$ \end{itemize}}

\newcommand{\ssnl}{\vskip 5pt} 
\newcommand{\nl}{\vskip 12pt} 

\newcommand{\rood}{\color{red}}
\newcommand{\blauw}{\color{blue}}
\newcommand{\zwart}{\color{black}}

\newcommand{\ot}{\otimes}

\newcommand{\tussenen}{\qquad\quad\text{and}\qquad\quad}

\newcommand{\inv}{^{-1}}

\numberwithin{thm}{section}   
\numberwithin{equation}{section} 

\newcommand{\mycomment}[2]{{\blauw {#1}\ssnl}{\zwart{#2}}}
\newcommand{\keepcomment}[1]{}
\newcommand{\oldcomment}[1]{}

\newenvironment{niets}{\begin{itemize}[leftmargin=15pt] \item[ ]\hspace{- 15 pt}\bf   \rm }{\end{itemize}}
\newcommand{\nul}{\begin{inspr}\begin{niets}}
\newcommand{\enul}{\end{niets}\end{inspr}}

\newcommand{\whl}{\widehat\Lambda}

\begin{document}



\centerline{\bf \Large Algebraic quantum groups and duality III} 
\vspace{5pt}
\centerline{\bf \large The modular and the analytic structure}
\vspace{13 pt}
\centerline{\it A.\ Van Daele \rm ($^*$)}
\vspace{20 pt}
{\bf Abstract} 
\nl
This is the last part of a series of three papers on the subject.
\ssnl
In the first part \cite{VD-part1} we have considered the duality of algebraic quantum groups. In that paper, we use the term \emph{algebraic quantum group} for a regular multiplier Hopf algebra with integrals. We treat the duality as studied in that theory. In the second part \cite{VD-part2} we have considered the duality for multiplier Hopf $^*$-algebras with positive integrals. The main purpose of that paper is to explain how it gives rise to a locally compact quantum group in the sense of Kustermans and Vaes.
\ssnl
In a preliminary section of \cite{VD-part2} we have mentioned the analytic structure of  such a $^*$-algebraic quantum group. We have not gone into the details, except for that part needed to obtain that the scaling constant is trivial and for proving that  the composition of a positive left integral with the antipode gives a right integral that is again positive. Also the construction of the square root 
of the modular element 
 is discussed and used to obtain the two Haar weights on the associated locally compact quantum group.
\ssnl
In this paper we treat the analytic structure of the $^*$-algebraic quantum group in detail. The analytic structure is intimately related with the modular structure. We obtain a collection of interesting formulas, relating the two structures. Further we compare these results with formulas obtained in the framework of general locally compact quantum groups.
\ssnl
As for the first two papers in this series, also here no new results are obtained. On the other hand  the approach is different, more direct and instructive. It might help the reader for understanding the more general theory of locally compact quantum groups.
\nl

Date: {\it 26 April 2023} 

\vskip 6cm
\hrule
\vskip 7 pt
\begin{itemize}
\item[($^*$)] Department of Mathematics, KU Leuven, Celestijnenlaan 200B,\newline
B-3001 Heverlee (Belgium). E-mail: \texttt{alfons.vandaele@kuleuven.be}
\end{itemize}

\newpage

\setcounter{section}{-1}  

\section{\hspace{-17pt}. Introduction} \label{sect:introduction}   

Let $(A,\Delta)$ be a multiplier Hopf $^*$-algebra with a positive right integral $\psi$. The space $A$ carries a scalar product, defined by  $\langle a',a\rangle =\psi(a^*a')$ where $a,a'$ are elements in $A$. One uses that the integral is automatically faithful. The completion of $A$ for this scalar product is denoted by $\mathcal H$ and $a\mapsto \Lambda(a)$ is used for the canonical embedding of $A$ in $\mathcal H$. It can be shown that there is a non-degenerate $^*$-representation $\pi$ of $A$ by \emph{bounded} operators on $\mathcal H$, given by $\pi(a)\Lambda(a')=\Lambda(aa')$ for all $a,a'\in A$. We think of this as the GNS-representation of $\psi$.
\ssnl
There is a unitary operator $V$ on $\mathcal H\ot \mathcal H$ satisfying $V(\Lambda(a)\ot\Lambda(a'))=(\Lambda\ot\Lambda)(\Delta(a)(1\ot a'))$  for $a,a'\in A$. One uses that $\Delta(A)(1\ot A)=A\ot A$ for a multiplier Hopf algebra $(A,\Delta)$ and that $\psi$ is right invariant. 
\ssnl
We now have the following result.

\stel
Let $M$  be the von Neumann algebra $\pi(A))''$. Define $\Delta_1(x)=V(x\ot 1)V^*$ for $x\in M$. Then $(M,\Delta_1)$ is a locally compact quantum group in the sense of Kustermans and Vaes. 
\estel

The coproduct on $M$ is the unique extension of the given coproduct on $A$. The left and right Haar weights on $M$ are obtained from the left and right integrals on $A$ using standard left Hilbert algebra techniques.
\ssnl
The result has been obtained in the literature and the procedure is reviewed in the second paper in this series of three, see \cite{VD-part2}, with some simplified and better methods.
\ssnl
Further in this introduction, we give more references.
\nl
In \cite{VD-part2}, we have mentioned the analytic structure of a $^*$-algebraic quantum group and proven some results, needed to construct the associated locally compact quantum group. In this third paper of the series, we give a new and more comprehensive treatment of the analytic structure with an emphasis on the modular properties of the invariant integrals. We find  the relations with the formulas for the dual objects that we reviewed in the first paper \cite{VD-part1}.  Moreover, we see how these formulas are translated into properties of the associated locally compact quantum group.
\ssnl
This \emph{series of three papers}, \cite{VD-part1,VD-part2} and this one, is meant to understand the theory of locally compact quantum groups, by the special case of a locally compact quantum group arising from a multiplier Hopf $^*$-algebra with positive integrals. The latter is of a purely algebraic nature and easier to understand. No new results are obtained, but the treatment is different, in some sense simpler, independent and still comprehensive.
\nl
\bf Content of these papers \rm
\nl
Before we give the content of this paper, we briefly also mention the content of the two preceding ones.
\ssnl
In the first paper \cite{VD-part1} we start with reviewing various properties of multiplier Hopf algebras $(A,\Delta)$ with integrals (algebraic quantum groups). We recall the definition of the associated objects, the left and right integral, the modular automorphisms, the modular element and the scaling constant. We collect the relations among these objects. In the second section of this first paper, we consider the duality for such an algebraic quantum group. One of the main issues is the relations of the objects of the dual with the ones of the original algebraic quantum group. Also the \emph{Fourier transform} is discussed here. In the third section, the duality is studied and seen as the canonical map $a\ot b\mapsto \Delta(a)(1\ot b)$. This is related with the action of the \emph{Heisenberg algebra}, generated by the algebra $A$ and its dual $B$, subject to the \emph{Heisenberg commutation relations}. \oldcomment{Mention the Heisenberg algebra.}
\ssnl
In the second paper \cite{VD-part2} we treat multiplier Hopf $^*$-algebras with positive integrals ($^*$-algebraic quantum groups). The main result is that a positive left integral exists if and only if there is a positive right integral. In fact, it is shown that the scaling constant has to be trivial and that $\psi\circ S$ is a positive left integral when $\psi$ is a positive right integral. The main purpose of that paper is the construction of the operator algebraic completions of the $^*$-algebraic quantum group. The difficulty lies in the construction of the Haar weights from the integrals and to show that these are still invariant. One uses a fair amount of the theory of left Hilbert algebras. Most of the results from this theory needed here are collected in two appendices of this second paper.
\ssnl
In the \emph{first section} of \emph{this paper}, we discuss the modular structure of a $^*$-algebraic quantum group. The main objects are the canonical involutions $\Lambda(a)\mapsto\Lambda(a^*)$ and $\whl(b)\mapsto \whl(b^*)$. Here $\Lambda$ is the canonical injection of $A$ into the GNS-Hilbert space $\mathcal H$ of the right integral $\psi$ on $A$ while $\whl$ is the canonical injection of the dual $B$ in $\mathcal H$ associated with the left integral $\widehat\varphi$ on $B$. 
\ssnl
In the \emph{second section}, we treat the analytic structure. The main property is that $A$ is spanned by eigenvectors, common for the maps $S^2$, $\sigma'$, $\sigma$, and left and right multiplication by $\delta$. Moreover all the eigenvalues are strictly positive. This allows to obtain one-parameter groups $\tau_t$, $\sigma'_t$ and $\sigma_t$ on the von Neumann completion of $A$ in the Hilbert space $\mathcal H$. Elements in $A$ are analytic and we have the expected analytic generators $S^2$, $\sigma'$ and $\sigma$. We also get the one-parameter groups $\delta^{it}$ of unitary elements in the multiplier algebra $M(A)$ and we have that $\delta$ is the analytic generator of this group. 
Recall that in \cite{VD-part2} we have already the square root $\delta^\frac12$. It is nothing else but the analytic extension of $t\mapsto \delta^{it}$ to the point $-\frac{i}{2}$.
\ssnl
Similar results are found for the dual $B$. 
\ssnl
Finally, in the \emph{appendix} we approach the results of Section \ref{sect:anal} from an other perspective. The duality of $A$ with $B$ gives rise to the right regular representation $V$ on the Hilbert space tensor product $\mathcal H\ot\mathcal H$. When $T$ and $\widehat T$ are the modular involutions, obtained as the closures of the maps  $\Lambda(a)\mapsto\Lambda(a^*)$ and $\whl(b)\mapsto \whl(b^*)$, then we have the relations $V^*(T\ot \widehat T)=(T\ot\widehat T)V$. Building further from this relation, we recover some of the properties obtained in the former sections in a different way.
\oldcomment{... to be completed}{}
\nl
\bf Notations and conventions \rm
\nl
We use the tensor product symbol in various cases. When we have a  $^*$-algebra $A$ (without any topology), then $A\ot A$ denotes the algebraic tensor product. It is again a $^*$-algebra in a natural way. When we have a von Neumann algebra $M$, we consider the von Neumann algebra tensor product $M\ot M$. Finally, $\mathcal H\ot \mathcal H$ denotes the Hilbert space tensor product if $\mathcal H$ is a Hilbert space. It should always be clear from the context which type we are using. If not, it will be explicitly stated. 
\ssnl

We use the {\it Sweedler notation} for a coproduct on an algebra, also in the case of a multiplier Hopf algebra $(A,\Delta)$ where it is not assumed that $\Delta$  maps $A$ into the tensor product $A\ot A$ but rather in its multiplier algebra $M(A\ot A)$. However, because we assume that $\Delta(a)(1\ot a')$ and $(a\ot 1)\Delta(a')$ are in $A\ot A$ for all $a,a'\in A$, the use of the Sweedler notation is still justified. 
\ssnl 
Remark however that the Sweedler notation is essentially only used as a means to write formulas and equations in a more readable way and not really as a tool to prove results. 
\oldcomment{For the same reason we will, in a few occasions, also use the Sweedler notation for coproducts on an operator algebra. However, in this case, we will just use it for explaining formulas, not in the proofs. Ik denk niet dat we dit doen.\rood Check!}{}
\ssnl
Finally, in order to avoid too many different notations, subscripts, etc.,  we use the same symbol for different objects. 

We believe this will not lead to confusions. In any case, whenever there is some possible doubt, we will be more explicit.

\ssnl
More information about notations and conventions are found in the corresponding items in the introductions of \cite{VD-part1, VD-part2}.
\nl
\bf Basic references \rm
\nl
For the theory of Hopf algebras, we refer to the well-known books by Abe \cite{Ab} and Sweedler \cite{Sw}. See also the more recent work by Radford \cite{R-bk}. The original work on multiplier Hopf algebras is \cite{VD-mha} and for multiplier Hopf algebras with integrals, it is \cite{VD-alg}. The use of the Sweedler notation for multiplier Hopf algebras has been explained in e.g.\ \cite{VD-tools} and more recently in \cite{VD-sw}.
\ssnl
Pairings of multiplier Hopf algebras have been first studied in \cite {Dr-VD}. Actions of multiplier Hopf algebras are studied in \cite{Dr-VD-Z}.
\ssnl
Locally compact quantum groups have been considered by various authors and, as we mentioned already in the preceding papers, there are also different approaches. Our main references here are \cite{K-V2} and \cite{K-V3}. See also \cite{VD-sigma} for a slightly simplified treatment of the theory. The reader is also advised to look at the more recent survey paper \cite{VD-warsaw}.
\ssnl
For the locally compact quantum groups arising from algebraic quantum groups, the first paper to consider is \cite{Ku-VD}. Related is the analytic structure of an algebraic quantum group, see \cite{Ku}.  A more recent treatment is found in \cite{DC-VD}. 
\nl
\bf Acknowledgments \rm
\nl
I am very grateful to my coauthor, M.B.\ Landstad and other colleagues and friends, both at the University of Trondheim and the University of Oslo (where part of these notes were developed) for the nice and fruitful atmosphere during my regular visits to these departments.

\section{\hspace{-17pt}. The modular structure of an algebraic quantum group} \label{sect:modu}

Let $(A,\Delta)$ be a multiplier Hopf $^*$-algebra with positive integrals, i.e.\ a $^*$-algebraic quantum group. Fix a positive right integral $\psi$. Consider its GNS-space. It is a Hilbert space $\mathcal H$ with an injective linear map $\Lambda:A\to \mathcal H$, with dense range and satisfying $\langle \Lambda(a),\Lambda(c)\rangle=\psi(c^*a)$ for all $a,c$ in $A$. There is a non-degenerate $^*$-representation $\pi$ of $A$ by \emph{bounded} operators, satisfying $\pi(a)\Lambda(x)=\Lambda(ax)$ for all $a,x\in A$. 
\ssnl
We use $(B,\Delta)$ for the dual $^*$-algebraic quantum group $(\widehat A,\widehat \Delta)$.
So we also use $\Delta$ for the coproduct on the dual. Similarly we use $S$ for the antipode and $\varepsilon$ for the counit, both for $A$ and for $B$. The other objects on the dual will systematically be denoted, covered with a $\ \widehat{}\,$-symbol. \oldcomment{Should we include this in the beginning of the previous section ?}
\ssnl
 For the associated pairing we also use $(a,b)\mapsto \langle a,b\rangle$. It will always be clear from the context when we use this notation for the scalar product on the Hilbert space $\mathcal H$ and when for the pairing of $A$ with $B$. 
\ssnl
There is a non-degenerate $^*$-representation $\gamma$ of $B$, again by \emph{bounded} operators on $\mathcal H$ given by $\gamma(b)\Lambda(x)=\sum_{(x)} \langle x_{(2)},b\rangle \Lambda(x_{(1)})$ when $b\in B$ and $x\in A$. 
\ssnl
We refer to \cite{VD-part2} for details about the above statements. Also recall the following definition.
\defin
Define $\widehat \varphi$ on $B$ by $\widehat\varphi(b)=\varepsilon(a)$ when $b=\psi(S(\,\cdot\,)a)$ for $a$ in $A$.
We think of $b$ as the Fourier transform $\widehat a$ of $a$.
\edefin

The linear functional $\widehat\varphi$ is a left integral on the dual $B$ of $A$. For all $a\in A$ we have $\widehat\varphi(\widehat a^*\widehat a)=\psi(a^*a)$. Hence $\widehat\varphi$ is also positive. The result is thought of as the Plancherel formula. It is used to realize the GNS-representation of $B$, associated with this positive linear functional  in the same Hilbert space $\mathcal H$. This is the content of the following proposition.

\prop
Given $b$ in $B$ we define $\whl(b)=\Lambda(a)$ when $b=\psi(S(\,\cdot\,)a)$ for $a\in A$. Then $\langle\whl(b),\whl(d)\rangle=\widehat\varphi(d^*b)$ for all $b,d\in B$ and $\gamma(b)\whl(y)=\whl(by)$ for all $b,y$ in $B$.
\eprop

All these results are found in \cite{VD-part1,VD-part2}. 
\ssnl
As we have done in \cite{VD-part2}, further in this paper, we will consider $A$ and $B$ as acting directly on $\mathcal H$. So we write $a\xi$ and $b\xi$ for $\pi(a)\xi$ and $\gamma(b)\xi$ whenever $a\in A$, $b\in B$ and $\xi\in \mathcal H$.
\nl
\bf The modular operators \rm
\nl
In \cite{VD-part2} we have the following results. We use $\sigma'$ for the modular automorphism of $\psi$ and $\widehat\sigma$ for the modular automorphism of $\widehat\varphi$. \oldcomment{Reference ?}{}

\prop
i) The map $\Lambda(a)\mapsto\Lambda(a^*)$ defined on $\Lambda(A)$ has a densely defined adjoint. Its closure is denoted by $T$. For all $a\in A$ the vector $\Lambda(a)$ belongs to the domain of the adjoint $T^*$ of $T$ and we have $T^*(\Lambda(a))=\Lambda(\sigma'(a^*))$.
\ssnl
ii) The map $\whl(b)\mapsto\whl(b^*)$ defined on $\whl(B)$ has a densely defined adjoint. Its closure is denoted by $\widehat T$.  For all $b\in B$ the vector $\whl(b)$ belongs to the domain of the adjoint $\widehat T^*$ of $\widehat T$ and we have $\widehat T^*\whl(b))=\whl(\widehat\sigma(b^*))$.
\eprop

\bew
i) Take $a,c\in A$. Then 
\begin{equation*}
\langle\Lambda(a^*),\Lambda(c)\rangle=\psi(c^*a^*)=\psi(a^*\sigma'(c^*))
\end{equation*}
and therefore $\langle\Lambda(a^*),\Lambda(c)\rangle=\langle \Lambda(\sigma'(c^*)),\Lambda(a)\rangle$. This is sufficient to prove the first statement.
\ssnl
ii) Similarly take $b,d\in B$. Then
\begin{equation*}
\langle\whl(b^*),\Lambda(d)\rangle=\widehat\varphi(d^*b^*)=\widehat\varphi(b^*\widehat\sigma(d^*))
\end{equation*}
and therefore $\langle\whl(b^*),\whl(d)\rangle=\langle \whl(\sigma(d^*)),\whl(b)\rangle$. This is sufficient for the second statement.
\ebew

Remark that these maps are conjugate linear and so we have
\begin{equation*}
 \langle T\xi,\eta\rangle = \langle \xi,T^*\eta\rangle^-=\langle T^*\eta,\xi \rangle
\end{equation*}
when $\xi$ belongs to the domain of $T$ and $\eta$ to the domain of $T^*$. We use $\lambda^-$ as a shorthand notation for the complex conjugate $\overline\lambda$ of a complex number $\lambda$.
\ssnl
We have $\whl(B)=\Lambda(A)$ so that $\widehat T$ and $\widehat T^*$ are also defined on $\Lambda(A)$. The formulas for the actions of the operators $\widehat T$ and $\widehat T^*$ on vector $\Lambda(a)$ with $a\in A$ are given in the following proposition. Here $\delta$ is the modular element of $(A,\Delta)$.

\prop
For all $a\in A$ we have
\begin{equation*}
\widehat T\Lambda(a)=\Lambda(S(a^*)\delta\inv)
\tussenen
\widehat T^*\Lambda(a)=\Lambda(S(a)^*).
\end{equation*}
\eprop

\bew
i) Let $b=\psi(S(\,\cdot\,)a)$ with $a\in A$. Then
\begin{align*}
\langle x, b^*\rangle
&=\langle S(x)^*,b\rangle^-\\
&=\psi(S(S(x)^*)a)^-=\psi(x^*a)^-\\
&=\psi(a^*x)=\psi(\delta a^*x\delta\inv) \\
&=\psi(S(\delta a^*x))=\psi(S(x)S(a^*)\delta\inv).
\end{align*}
We have used that $\sigma'(\delta)=\delta$ because the scaling constant is $1$. We also have used $\psi(S(a))=\psi(a\delta\inv)$ and finally also that $S(\delta)=\delta\inv$. We get that $b^*=\psi(S(\,\cdot\,)S(a^*)\delta\inv)$. Then we get $\widehat T\Lambda(a)=\Lambda(S(a^*)\delta\inv)$.
\ssnl
ii) We now calculate the adjoint of this map. We have
\begin{align*}
\langle \Lambda(x),\Lambda(S(y^*)\delta\inv)\rangle
&=\psi(\delta\inv S\inv(y)x))\\
&=\psi(S(\delta\inv S\inv(y)x)\delta\inv)\\
&=\psi(S(x)y)=\langle \Lambda(y),\Lambda(S(x)^*\rangle.
\end{align*}
This means that 
\begin{equation*}
\langle\Lambda(x),\widehat T\Lambda(y)\rangle=\langle\Lambda(y),\Lambda(S(x)^*)\rangle
\end{equation*}
and therefore $\widehat T^*\Lambda(x)=\Lambda(S(x)^*)$.
\ebew

For completeness, we also calculate $T\widehat\Lambda(b)$ and $T^*\widehat\Lambda(b)$. 

\prop\label{prop:1.5}
 For $b\in B$ we have $T\whl(b)=\whl(S(b)^*\widehat \delta)$ and $T^*\whl(b)=\whl(S(b^*))$.
\eprop
\bew
i) As before we let $b=\psi(S(\,\cdot\,)a)$. Then
\begin{equation*}
T\whl(b)=T\Lambda(a)=\Lambda(a^*)=\whl(d)
\end{equation*}
where $d=\psi(S(\,\cdot\,)a^*)$. So we get for all $x\in A$
\begin{align*}
\langle x, d\rangle
&=\psi(S(x)a^*) = \psi(aS(x)^*)^- \\
&=\psi({\sigma'}\inv(S(x)^*)a)^-\\
&=\langle {\sigma'}\inv (S(x)^*),S\inv b\rangle^-\\
&=\langle S(x)^*,\widehat\delta S^2(S\inv(b))\rangle^-\\
&=\langle S(x)^*, \widehat \delta S(b)\rangle^-\\
&=\langle x, S(b)^*\widehat\delta\rangle.
\end{align*}
We have used a formula from Proposition 2.5 of \cite{VD-part1}. We see that $d=S(b)^*\widehat\delta$.
\ssnl
ii)
Again with $b=\psi(S(\,\cdot\,)a)$ we find
\begin{equation*}
T^*\whl(b)=T^*\Lambda(a)=\Lambda(\sigma'(a^*))=\whl(d)
\end{equation*}
where now $d=\psi(S(\,\cdot\,)\sigma'(a^*))$. Then for all $x\in A$
\begin{align*}
\langle x,d\rangle
&=\psi(S(x)\sigma'(a^*))=\psi(a^*S(x))\\
&=\psi(S(x)^*a)^-\\
&=\langle S(x)^*,S\inv(b)\rangle^-\\
&=\langle x,S(b^*)\rangle
\end{align*}
and we find $d=S(b^*)$.
\ebew
\oldcomment{Check!}{}
\defin
We consider the polar decompositions of these maps
\begin{equation*}
T=J\nabla^\frac12 \tussenen \widehat T=\widehat J \widehat\nabla^\frac12.
\end{equation*}
\edefin

We have the following realizations.

\prop\label{prop:1.6}
\begin{align*}
\nabla\Lambda(a))&=\Lambda(\sigma'(a))
&\nabla\widehat\Lambda(b)=\widehat\Lambda(S^2(b)\widehat\delta\inv)\\
\widehat\nabla\widehat\Lambda(b))&=\widehat\Lambda(\widehat\sigma(b))
&\widehat\nabla\Lambda(a)=\Lambda(S^{-2}(a)\delta).
\end{align*}
\eprop
\bew
i) For $a\in A$ we have 
\begin{equation*}
\nabla\Lambda(a)=T^*T\Lambda(a)=T^*\Lambda(a^*)=\Lambda(\sigma'(a)).
\end{equation*}
Similarly for $b\in B$ we get
\begin{equation*}
\widehat\nabla\whl(b)=\widehat T^*\widehat T\whl(b)=\widehat T^*\whl(b^*)=\whl(\widehat\sigma(b)).
\end{equation*}
\ssnl
ii) When $a\in A$ we get
\begin{align*}
\widehat\nabla(\Lambda(a))
&=\widehat T^*\Lambda(S(a^*)\delta\inv)\\
&=\Lambda(S(S(a^*))^*S(\delta\inv)^*)\\
&=\Lambda(S^{-2}(a)\delta).
\end{align*}
Finally we calculate $\nabla\whl(b)$ for $b\in B$.
Given $b=\psi(S\,\cdot\,)a)$ we have
\begin{equation*}
\nabla\widehat\Lambda(b)=\nabla\Lambda(a)=\Lambda(\sigma'(a))=\Lambda(d)
\end{equation*}
where $d=\psi(S(\,\cdot\,)\sigma'(a))$.
We find
\begin{align*}
\langle x, d\rangle
&=\psi(S(x)\sigma'(a))=\psi({\sigma'}\inv(S(x))a)\\
&=\psi(S(\sigma(x))a)= \langle\sigma(x),b\rangle.
\end{align*}
Now we have$\langle\sigma(x),b\rangle=\langle x,S^2(b)\widehat \delta\inv\rangle$. This is the dual version of a formula we have in Proposition 2.5 of \cite{VD-part1}. It follows that $\nabla\widehat\Lambda(b)=\widehat \Lambda(S^2(b)\widehat\delta\inv)$.
\ebew

The last formula can be verified using the formulas in Proposition \ref{prop:1.5}. Indeed
\begin{equation*}
T^*T\whl(b)=T^*\whl(S(b)^*\widehat \delta)=\whl(S((S(b^*)\widehat\delta)^*)=\whl(S^2(b){\widehat \delta}\inv)
\end{equation*}

We obviously can check the following implementations.

\prop\label{prop:1.7}
\begin{align*}
\nabla a\nabla\inv&=\sigma'(a)
&\nabla b\nabla \inv = S^2(b)\mbox{\ \  } \\
\widehat\nabla b\widehat\nabla\inv& =\widehat\sigma(b)
&\widehat\nabla a \widehat\nabla\inv=S^{-2}(a)
\end{align*}
\eprop

In the appendix, see Proposition \ref{prop:JIV} we have \emph{one-parameter} forms of these formulas. 
\ssnl
For the action of $\delta$ and $\widehat \delta$ we find the following formulas.

\prop
\begin{equation*}
\widehat\delta\Lambda(a)=\Lambda(S^2\sigma\inv(a))
\tussenen
\delta\whl(b)=\whl(S^2\widehat\sigma'(b))
\end{equation*}
\eprop
\bew
i) For the first formula, we use the definition of the action of $B$ on $\Lambda(a)$. Then we find
\begin{align*}
\widehat\delta\Lambda(a)
&=\sum_{(a)}\langle a_{(2)},\widehat\delta\rangle \Lambda(a_{(1)})\\
&=\sum_{(a)}\varepsilon(\sigma\inv( a_{(2)}) \Lambda(a_{(2)})\\
&=\sum_{(a)}\varepsilon(\sigma\inv( a_{(2)}) \Lambda(S^2S^{-2}(a_{(2)}))\\
&=\sum_{(\sigma\inv(a))}\varepsilon((\sigma\inv(a)_{(2)}) \Lambda(S^2(\sigma\inv(a)_{(2)}))\\
&=\Lambda(S^2\sigma\inv(a)).
\end{align*}
We have used that $\Delta(\sigma(a))=(S^2\ot\sigma)\Delta(a)$ (Proposition 1.13 in \cite{VD-part1}) and $\langle a,\widehat \delta\rangle=\varepsilon(\sigma\inv(a))$ (Proposition 2.3 in \cite{VD-part1}) for all $a$.
\ssnl
ii) As before let $b=\psi(S(\,\cdot\,)a)$ with $a\in A$. Then $\delta\whl(b)=\delta\Lambda(a)=\Lambda(\delta a)$ so that $\delta\whl(b)=\whl(d)$ where $d=\psi(S(\,\cdot\,)\delta a)$. For all $x\in A$ we  have
\begin{align*}
\langle x,d\rangle
&= \psi(S(x)\delta a)=\langle S(x)\delta, S\inv(b)\rangle\\
&=\langle S^{-2}(S^3(x))\delta,S\inv(b)\rangle=\langle S^3(x),\widehat\sigma\inv S\inv(b)\rangle\\
&=\langle x, S^3\widehat\sigma\inv S\inv(b)\rangle=\langle x,S^2\widehat\sigma'(b)\langle.
\end{align*}
Here we have used that $\langle a,\widehat\sigma(b)\rangle=\langle S^2(a)\delta\inv(a),b\rangle$ for all $a,b$ (Proposition 2.5 in \cite{VD-part1}) and $S\widehat\sigma\inv(a)=\widehat\sigma'(S(a))$ for all $a$ (Proposition 1.10 in \cite{VD-part1}).
\ebew

We finish this section with the following remark. 

\opm
As we have mentioned, there is a problem when comparing the formulas above with those found in the theory of locally compact quantum groups. The reason for that is that twofold. On the one hand, there is the different convention for the coproduct on the dual. On the other hand, we start with the right integral and not with the left integral. Considering the latter, we might replace $\Delta$ on $A$ by $\Delta^\text{op}$. The effect is that $S$ on $A$ is replaced by $S\inv$,  that $\sigma'$ is replaced by $\sigma$ and $\delta$ by $\delta\inv$. This further leads to a replacement of the product, as well as of the coproduct on $B$. 
\eopm

\section{\hspace{-17pt}. The analytic structur of an algebraic quantum group} \label{sect:anal} 

As in the previous section, we have a $^*$-algebraic quantum group $(A,\Delta)$  with a positive right integral $\psi$. A left integral on $A$ is given by $\varphi=\psi\circ S$. The modular automorphisms of $\psi$ and $\varphi$ are denoted by respectively $\sigma'$ and $\sigma$. 
\nl
In \cite{VD-part2} we have shown the following result, see Proposition 1.6 in \cite{VD-part2}. \oldcomment{Check reference!}{}

\prop\label{prop:1.1}
 For all $a$ in $A$ there is is finite-dimensional subspace containing all element $\delta^na$ and $a\delta^n$ for all integers $n$.
\eprop 

An easy consequence of the formulas relating the various objects gives that the maps $S^2, \sigma,\sigma'$ mutually commute (see Proposition 1.10 in \cite{VD-part1}). 
Denote $S^{-2}\sigma$  by $\kappa$ and $S^2\sigma'$ by $\rho$. They are still isomorphisms (but in general not $^*$-isomorphisms).
\ssnl
In Proposition 1.12 of \cite{VD-part1} we have the formulas 
\begin{equation*}
\Delta(S^2(a))=(S^2\ot S^2)\Delta(a)
\tussenen
\Delta(S^2(a))=(\sigma\ot{\sigma'}\inv)\Delta(a)
\end{equation*}
for all $a$. Then it follows that $(\kappa\ot\rho\inv)\Delta(a)=\Delta(a)$ for all $a$. We use this to prove the following result, see Proposition 2.3 in \cite{DC-VD}.
\oldcomment{\ssnl We need to refer to the work with Kenny }{}

\prop
For all $a\in A$ there exists a finite-dimensional subspace of $A$ containing all the elements $\kappa^n(a)$ and $\rho^n(a)$ for all integers $n$.
\eprop

\bew
i) Take $a,b,c$ in $A$ and write $a\ot b=\sum_i \Delta(p_i)(1\ot q_i)$. Apply $\kappa^n\ot \rho^{-n}$ and multiply with $c$ to get
\begin{equation*}
\kappa^n(a)\ot c\rho^{-n}(b)=\sum_i(1\ot c)\Delta(p_i)(1\ot \rho^{-n}(q_i).
\end{equation*}
Using that $(1\ot c)\Delta(p_i)$ belongs to $A\ot A$ we have a finite-dimensional subspace $L$ of $A$ such that $(1\ot c)\Delta(p_i)\in L\ot A$ for all $i$. We see that $\kappa^n(a)$ belongs to $L$. 
\ssnl
Similarly, if we start with $a\ot b=\sum_i \Delta(p_i)(q_i\ot 1)$ we get a finite-dimensional subspace of $A$ containing all elements $\rho^n(a)$.
\ebew

We can also give the following argument. Using that 
\begin{equation*}
\Delta(\sigma(a))=(S^2\ot\sigma)\Delta(a)
\tussenen
\Delta(S^2(a))=(S^2\ot S^2)\Delta(a)
\end{equation*}
we find that $\Delta(\kappa)=(\iota\ot\kappa)\Delta(a)$. Then we write $\Delta(a)(1\ot c)=\sum_i p_i\ot q_i$, apply $\iota\ot \kappa^n$ and $\varepsilon$ to arrive at
\begin{equation*}
\kappa^n(a)\varepsilon(\kappa^n(c))=\sum_i p_i\varepsilon(\kappa^n(q_i).
\end{equation*}
Now we have $\varepsilon((\kappa^n(b))=\langle b,{\widehat\delta}\inv\rangle$ and by choosing $c=\widehat\varphi(\,\cdot\,b^*b)$ we get 
\begin{equation*}
\langle b,{\widehat\delta}\inv\rangle=\widehat\varphi({\widehat\delta}^nb^*b).
\end{equation*}
This will be non-zero as soon as $b$ is non-zero. 
In the end we get that $\kappa^n(a)$ belongs to the space spanned by the elements $p_i$. A similar argument works for $\rho$.
\keepcomment{\ssnl Is this useful? Should we include this or just keep it for our own information?}{}
\ssnl
We now combine all these results.

\prop
For all $a$ there is a finite-dimensional subspace of $A$ containing all elements $S^n(a)$, $\sigma^n(a)$ and ${\sigma'}^n(a)$ for all integers $n$. 
\eprop
\bew
We have $\sigma\sigma'=\kappa\rho$ because $\kappa=S^{-2}\sigma$ and $\rho=S^2\sigma'$. This will give that $\sigma^n{\sigma'}^n(a)$ belongs to a finite-dimensional space. Further we have $\sigma\sigma'(a)=\delta\sigma^2(a)\delta\inv$. This, combined with Proposition \ref{prop:1.1} will give that $\sigma^n(a)$ belongs to a finite-dimensional space. Then a combination of all these properties will  complete the proof of this proposition.
\ebew

We have these maps on these finite-dimensional spaces. The following will show that they have positive eigenvalues.
\prop
For all $a$ we have
\begin{equation*}
\psi(a^*S^2(a))\geq 0,
\qquad
\psi(a^*\sigma(a))\geq 0,
\qquad
\psi(a^*\sigma'(a))\geq 0.
\end{equation*}
\eprop
\bew
i) For all $a$ we have
\begin{equation*}
\psi(a^*S^2(a))=\varphi(S(a)S\inv(a^*))=\varphi(S(a)S(a)^*).
\end{equation*}
Because $\varphi$ is also positive, see Theorem 1.7 in \cite{VD-part2}, we find that $\psi(a^*S^2(a))\geq 0$.
\ssnl
ii) For all $a$ we have
\begin{equation*}
\psi(a^*\sigma'(a))=\psi(aa^*)\geq 0.
\end{equation*}
iii) Finally, again for all $a$ we have, because $\sigma(\delta)=\delta$, 
\begin{equation*}
\psi(a^*\sigma(a))=\varphi(a^*\sigma(a)\delta)=\varphi(a\delta a^*)\geq 0.
\end{equation*}

\vskip -10pt\ebew

Now we can obtain the main result that will eventually yield the analytic structure and the results for the completion. See Theorem 3.5 in \cite{DC-VD}. \oldcomment{Include also multiplication with $\delta$.}{}

\prop
The space $A$ is spanned by common eigenvectors of the maps $S^2$, $\sigma$ and $\sigma'$.
\eprop

\bew
i) Let $a\in A$ and consider the space $L$ spanned by the elements $S^2(a)$. It is finite-dimensional and invariant under $S^2$. For the scalar product induced by $\psi$, the map $S^2$ is self-adjoint and positive. It follows that $L$ is spanned by eigenvectors of $S^*$ and that the eigenvectors are all strictly positive. As this is true for all $a$ we get that $A$ is spanned by eigenvectors of $S^2$. In a similar way, $A$ is spanned by eigenvectors for $\sigma$ and eigenvectors for $\sigma'$. 
\ssnl
ii) Because all these maps commute, we can get a basis of vectors that are eigenvectors jointly for the three maps $S^2$ $\sigma$ and $\sigma'$.
\ebew

We have seen already that $A$ is also spanned by eigenvectors of left and right multiplication by $\delta$. This was proven already in Proposition 1.7 of \cite{VD-part2}. Because we  have $\psi(a^*\delta a)=\psi((\delta^\frac12 a)^*\delta^\frac12 a)$ we see that also $\psi(a^*\delta a)\geq 0$. On the other hand $\varphi(a^*a\delta))=\psi(a^*a)\geq 0$. Therefore also the eigenvalues of multiplication with $\delta$ are strictly positive. Because these maps also commute with the maps $S^2$,  $\sigma$ and 
$\sigma$ and $\sigma'$, we have even the following result.

\prop
The space $A$ is spanned by common eigenvectors of the maps $S^2$, $\sigma$ and $\sigma'$ and left and right multiplication by $\delta$. All eigenvalues are strictly positive.
\eprop 

\oldcomment{To be  completed with reference to the original papers and with my earlier papers. Perhaps also formulate better.}{} 

By duality, we also have this result for the dual.

\prop
The space $B$ is spanned by common eigenvectors of the maps $S^2$, $\widehat\sigma$ and $\widehat\sigma'$ and left and right multiplication by $\widehat\delta$. All eigenvalues are strictly positive.
\eprop

Again see Theorem 3.5 of \cite{DC-VD}.
\nl
\bf The one-parameter extensions \rm
\nl
We can now define various one-parameter groups. First we have the powers of $\delta$. 

\prop
We can define multiplier $\delta^{it}$ for all $t\in\mathbb R$ by 
\begin{equation*}
\delta^{it}a=\lambda^{it}a
\tussenen
c\delta^{it}=\mu^{it}c
\end{equation*}
where $a$ and $c$ are elements in $A$ satisfying 
$\delta a=\lambda a$ and $c\delta=\mu c$ for some strictly positive numbers $\lambda$ and $\mu$.
\eprop

\bew
We proceed as in the proof of Proposition 1.11 of \cite{VD-part2}.
\ssnl
Let $t\in \mathbb R$. Take $a$ in $A$ and write it as a sum $\sum_i a_i$ where $\delta a_i=\lambda_i a_i$. We assume that all $\lambda_i$ are different. Then the decomposition is unique because different  elements $a_i$ now satisfy $\psi(a_i^*a_j)=0$ when $i\neq j$. Indeed
\begin{equation*}
\psi(a_i^*\delta a_j)=\lambda_j \psi(a_i^*a_j)
\tussenen 
\psi((\delta a_i)^*a_j)=\lambda_i\psi(a_i^*a_j)
\end{equation*}
and $\psi(a_i^*\delta a_j)=\psi((\delta a_i)^*a_j)$.
Then we can define $\delta^{it}a:=\sum_i \lambda^{it}a_i$. 
\ssnl
In a similar way we can define $a\delta^{it}$. It is true that $(\delta^{it}a)\delta^{it}=\delta^{it}(a\delta^{it})$ so that we do have a well-defined multiplier of $A$.
\ebew

They are unitary operators and can defined on all of the Hilbert space $\mathcal H$.
\ssnl
We have analytic extensions and $\delta^n=\delta^{iz}$ when $z=-in$. 
\ssnl
\oldcomment{We should write a proof in detail for ourselves - Just to be sure!}{}

\prop 
We also have $\Delta(\delta^{it})=\delta^{it}\ot \delta^{it}$ for all $t$.
\eprop
\bew
Take elements $a,c$ satisfying $\delta a=\lambda a$ and $\delta c=\mu c$. Write $a\ot c=\sum_i \Delta(p_i)(1\ot q_i)$. Then multiply from the left with $\delta\ot\delta$ and use that $\Delta(\delta)=\delta\ot\delta$. We get $\lambda\mu (a\ot c)=\sum_i \Delta(\delta p_i)(1\ot q_i)$ This implies that $\delta p_i=\lambda\mu p_i$ for all $i$. Then 
\begin{align*}
(\delta^{it}\ot \delta^{it})(a\ot c) 
&= \lambda^{it}\mu^{it}(a\ot c) \\
&= \sum_i \Delta((\lambda\mu)^{it}p_i)(1\ot q_i)\\
&= \sum_i \Delta(\delta^{it}p_i)(1\ot q_i) \\
&=\Delta(\delta^{it})\sum_i\Delta(p_i)(1\ot q_i).
\end{align*}

\ebew

For the one-parameter groups of modular automorphisms, we will use another approach.
\ssnl
Recall that we have the self-adjoint operators $\nabla$ and $\widehat\nabla$ on the Hilbert space $\mathcal H$. They give rise to unitary one-parameter groups of unitaries $\nabla^{it}$ and $\widehat\nabla^{it}$. The space $\mathcal H$ is the closed linear span of eigenvectors, all contained in $\Lambda(A)$. We have $\nabla^{it}\Lambda(a)=\lambda^{it}\Lambda(a)$ for elements $a$ satisfying $\sigma'(a)=\lambda a$ and $\widehat\nabla^{it}\whl(b)=\mu^{it}b$ for elements $b$ satisfying $\widehat\sigma(b)=\mu b$.

\prop
\begin{equation*}
\nabla^{it} A \nabla^{-it} \subseteq A
\tussenen
\widehat\nabla^{it} B \widehat\nabla^{-it} \subseteq B.
\end{equation*}
\eprop
\bew
We use the formulas obtained in Proposition \ref{prop:1.6}.
\ssnl
i) Assume that $a,c\in A$ and that $\sigma'(a)=\lambda a$ and $\sigma'(c)=\mu c$. Then $\sigma'(ac)=\lambda\mu ac$  and we find
\begin{align*}
\nabla^{it}a\Lambda(c)
&=\nabla^{it}\Lambda(ac) \\
&=\Lambda(\lambda^{it}\mu^{it}(ac)) \\
&=\lambda^{it}\nabla^{it}\Lambda(c).
\end{align*}
This holds for all such elements $c$ and hence for all $c\in A$. It follows that $\nabla^{it}a\nabla^{-it}=\lambda^{it}a$. Because $A$ is spanned by eigenvectors, we get $\nabla^{it}A\nabla^{-it}\subseteq A$.
\ssnl
ii) In a similar way we get $\widehat\nabla^{it} B \widehat\nabla^{-it} \subseteq B$.
\ebew

Then we have the following definitions.

\defin
We can define $\sigma'_t$ on $A$ for all $t$ by $\sigma'_t(a)=\nabla^{it}a\nabla^{-it}$ and $\widehat \sigma_t$ on $B$ for all $t$  by $\widehat \sigma_t(b)=\widehat\nabla^{it} b \widehat\nabla^{-it}$. 
\edefin

It is also obvious to show the following.

\prop
The elements of $A$ are analytic for the one-parameter group $\sigma'_t$ and $\sigma'_{-i}(a)=\sigma'(a)$ for $a\in A$. Similarly, elements of $B$ are analytic for the one-parameter group $\widehat\sigma_t$ and $\widehat\sigma_{-i}(b)=\widehat\sigma(b)$ for $b\in B$.
\eprop

Similarly we have the following.

\prop\label{prop:2.13}
\begin{equation*}
\widehat\nabla^{it} A \widehat\nabla^{-it} \subseteq A
\tussenen
\nabla^{it} B \nabla^{-it} \subseteq B.
\end{equation*}
\eprop
\bew
We use the formulas obtained in Proposition \ref{prop:1.6}. 
 \ssnl
 i) Let $a$ be and element in $A$ such that $S^{-2}(a)=\lambda a$, Further let $c\in A$ satisfying $S^{-2}(c)\delta=\mu c$. Then we get
 \begin{equation*}
\widehat\nabla(\Lambda(ac))=\Lambda(S^{-2}(ac)\delta)=\lambda\mu \Lambda(ac).
\end{equation*}
It follows that 
 \begin{equation*}
\widehat\nabla^{it}(\Lambda(ac))=\lambda^{it}\mu^{it} a \Lambda(c)= \lambda^{it}\Lambda(S^{-2}(c)\delta)=\lambda^{it}\widehat\nabla^{it}\Lambda(c).
\end{equation*}
Then this holds for all $c\in A$ and we get
\begin{equation*}
\widehat\nabla^{it} a \widehat\nabla^{-it} =\lambda^{it}a.
\end{equation*}
Then the result follows because $A$ is spanned by such eigenvectors.
\ssnl
ii) The other statement follows in a similar way from the formula $\nabla\whl(b)=\whl(S^2(b)\delta\inv)$.

\ebew

\defin
We now define $\tau_t$ on $A$ for all $t$ by $\tau_t(a)=\widehat\nabla^{it} a \widehat\nabla^{-it}$ and we define $\widehat\tau_t$ on $B$ for all $t$ by $\widehat\tau_t(b)=\nabla^{it}b\nabla^{-it}$.
\edefin

Again we have that elements of $A$ and $B$ are analytic and $\tau_{-i}(a)=S^{-2}(a)$ while $\widehat\tau_{-i}(b)=S^2(b)$. Remark that different conventions are used, see the discussion in \cite{VD-warsaw} about this problem.

\prop
For these one-parameter groups, we have for all $x\in M$ and $y\in \widehat M$, 
\begin{align*}
\Delta(\sigma'_t(x))&=(\sigma'_t\ot \tau_t)\Delta(x)
&\Delta(\widehat\sigma_t(y))=(\widehat\tau_t\ot \widehat\sigma_t)\Delta(y)\\
\Delta(\tau_t(x))&=(\tau_t\ot\tau_t)\Delta(x)
&\Delta(\widehat\tau_t(y))=(\widehat\tau_t\ot\widehat\tau_t)\Delta(y)
\end{align*}
\eprop

\section{\hspace{-17pt}. Conclusions and further remarks} \label{sect:concl} 

In a series of three paper, \cite{VD-part1},\cite{VD-part2} and this one, we have given an overview of the theory of algebraic quantum groups and their duality. 
\ssnl
In the first one, we have collected and discussed the formulas for a multiplier Hopf algebra $(A,\Delta)$ and its dual $(\widehat A,\widehat \Delta)$. In the second paper, we have studied multiplier Hopf $^*$-algebras with positive integrals and we have shown how it can be completed to a locally compact quantum group. Finally, in this paper, the focus lied on the modular properties of the integrals on $A$ and on $B$, together with the analytic properties. 
\ssnl
No new results are obtained, but the treatment is newer, sometimes very different, more independent and still complete. 
\ssnl
The purpose of writing these papers is for those who want to learn about the general theory of locally compact quantum groups, but starting with the simpler case of a locally compact quantum group obtained from a $^*$-algebraic quantum group.
\ssnl
\oldcomment{To be completed? Something about new research?}{}
\oldcomment{Relation with our work on bicrossproduct, it was the motivation for writing these papers.}{}




\renewcommand{\thesection}{\Alph{section}}

\setcounter{section}{0}



\renewenvironment{stelling}{\begin{itemize}\item[ ]\hspace{-28pt}\bf Theorem \rm }{\end{itemize}}
\renewenvironment{propositie}{\begin{itemize}\item[ ]\hspace{-28pt}\bf Proposition \rm }{\end{itemize}}
\renewenvironment{lemma}{\begin{itemize}\item[ ]\hspace{-28pt}\bf Lemma \rm }{\end{itemize}}

\section{\hspace{-17pt}. Appendix. A multiplicative unitary approach}\label{sect:appA}   

We have a multiplier Hopf $^*$-algebra $(A,\Delta)$ with a positive right integral $\psi$. We consider the associated left integral $\widehat\varphi$ on the dual $(B,\Delta)$ defined by $\widehat\omega(b)=\varepsilon(a)$ when $b=\psi(S(\,\cdot\,)a)$ for $a\in A$. The GNS-space of $\psi$ is $\mathcal H$ with the canonical embedding $\Lambda$ satisfying $\langle\Lambda(a),\Lambda(c)\rangle=\psi(c^*a)$. Also the GNS-space of $\widehat\varphi$ is $\mathcal H$ and the canonical embedding $\whl$ of $B$ in $\mathcal H$ is given by $\whl(b)=\Lambda(a)$ when again $b=\psi(S(\,\cdot\,)a)$ for $a\in A$. This is possible because we have the Plancherel formula saying that $\widehat\varphi(b^*b)=\psi(a^*a)$ when $b=\psi(S(\,\cdot\,)a)$.
\ssnl
Recall the following result. 
\prop
There is a unitary operator $V$ on $\mathcal H\ot \mathcal H$ defined by
\begin{equation*}
V(\Lambda(x)\ot\xi)=\sum_{(x)} \Lambda(x_{(1)})\ot x_{(2)}\xi.
\end{equation*}
when $x\in A$ and $\xi\in\Lambda(A)$. It satisfies 
\begin{equation*}
V^*(\xi\ot\widehat\Lambda(y)=\sum_{(y)}y_{(1)}\xi\ot\widehat\Lambda(y_{(2)})
\end{equation*}
when $y\in B$ and $\xi\in\widehat\Lambda(B)$.
\eprop

The first result is shown in Proposition 2.5 of \cite{VD-part2} while the second one is found in Proposition 4.7 of that paper. 
\ssnl
Remark that 
$V$ belongs to $M(\widehat A\ot A)$. The right leg generates $A$ and the von Neumann algebra $M$ while the left leg generates $\widehat A$ and the von Neuman algebra $\widehat M$. Then $V$ sit in the von Neumann tensor product $\widehat M\ot M$. 
\ssnl
We have defined $T$ and $\widehat T$ as the closures of the maps $\Lambda(a)\to\Lambda(a^*)$ and $\whl(b)\mapsto \whl(b^*)$ with $a\in A$ and $b\in B$. We have used
\begin{equation*}
T=J\nabla^\frac12
\tussenen
\widehat J\widehat \nabla^\frac12
\end{equation*}
for the polar decompositions of these maps.
\ssnl
In Proposition 3.1 of  \cite{VD-part2} we have proven the following result. 

\prop
For all $a$ we have that $\Lambda(a)$ belongs to the domain of $T^*$ and
 $T^*(\Lambda(a)= \Lambda(\sigma'(a^*))$ for all $a\in A$..
\eprop

In Proposition 4.9 of \cite{VD-part2} we also consider the dual version of this, but now we obtain the following property.

\prop
We have $\widehat T(\Lambda(a))=\Lambda(S(a^*)\delta\inv))$ for all $a\in A$. Also $\Lambda(a)$ belongs to the domain of $\widehat T^*$ and  $\widehat T^*\Lambda(a)=\Lambda(S(a)^*)$.
\eprop 

Using these to formulas, we get the following relation of these operators with $V$.

\prop
We have $V^*(T\ot \widehat T)=(T\ot \widehat T)V$.
\eprop
\bew
For all $x,y\in A$ we get
\begin{align*}
V^*(T\ot \widehat T)(\Lambda(x)\ot\Lambda(y))
&=V^*(\Lambda(x^*)\ot\Lambda(S(y^*)\delta\inv)\\
&=\sum_{(x)} \Lambda(x_{(1)}^*)\ot S(x_{(2)}^*)\Lambda(S(y^*)\delta\inv)\\
&=\sum_{(x)} \Lambda(x_{(1)}^*)\ot \Lambda(S(x_{(2)}^*)S(y^*)\delta\inv)\\
&=\sum_{(x)} \Lambda(x_{(1)}^*)\ot \Lambda(S((x_{(2)}y)^*)\delta\inv)\\
&=\sum_{(x)} (T\ot \widehat T)(\Lambda(x_{(1)})\ot \Lambda(x_{(2)}y))\\
&=(T\ot \widehat T)V(\Lambda(x)\ot\Lambda(y)).
\end{align*}
\ebew

\opm
In earlier work we used the map $\Lambda(a)\mapsto\Lambda(S(a^*))$ instead of $\widehat T$. This also implements $a\mapsto S(a^*)$ and that is what we need for this formula. With $\widehat T$ we make a better choice.
\eopm

From the uniqueness of the polar decompositions we get:

\prop\label{prop:JIV}
We have 
\begin{equation}
(J\ot \widehat J)V=V^*(J\ot \widehat J)
\tussenen
(\nabla^{it}\ot {\widehat\nabla}^{it})V=V(\nabla^{it}\ot {\widehat\nabla}^{it}).
\end{equation}
\eprop

In the first place we have the modular data. For $x\in M$ and $y\in \widehat M$ we have
\begin{align*}
\sigma'_t(x)&=\nabla^{it}x\nabla^{-it}\in M \tussenen x\mapsto Jx^*J\in M' \\
\widehat\sigma_t(y)&=\widehat\nabla^{it}y\widehat\nabla^{-it}\in\widehat M \tussenen y\mapsto \widehat Jy^*\widehat J\in \widehat M.
\end{align*}

\prop
For $x\in M$  and $y\in \widehat M$ we define
\begin{align*}
\tau_t(x)&=\widehat\nabla^{-it}x\widehat\nabla^{it}\in M
\tussenen R(x)=\widehat Jx^*\widehat J\in M\\
\widehat\tau_t(y)&=\nabla^{it}y\nabla^{-it}\in \widehat M
\tussenen \widehat R(y)=Jy^*J\in \widehat M.
\end{align*}
Then $R$ is an involutive anti-isomorphism of $M$ while $\widehat R$ is an involutive anti-isomorphism of $\widehat M$. The $\tau$-maps are one parameter groups of $^*$-automorphisms.
\eprop

\bew
Because the right leg of $V$ is dense in $M$, it follows from the first formula in Proposition \ref{prop:JIV} that $R$ maps $M$ to itself. 
Similarly for $\widehat R$ because the left leg of $V$ is dense in $\widehat M$.
\ebew

We call $R$ and $\widehat R$ the unitary antipode on $M$ and $\widehat M$ respectively. The automorphisms $\tau$ and $\widehat\tau$ are the scaling groups.
\oldcomment{We need a comment on the choices we make for $\tau_t$ and $\widehat\tau_t$, and refer to the different conventions.}{}
\nl
\bf Properties of these maps \rm
\nl
First we use that $V$ implements the coproducts. This gives the following result. \mycomment{Aan te passen aan de conventies!}{}

\prop
For $x\in M$ and $y\in\widehat M$ we have
\begin{align*}
\Delta(\sigma'_t(x))&=(\sigma'_t\ot \tau_{-t})\Delta(x)\\
\Delta(\widehat\sigma_t(y))&=(\widehat\tau_t\ot \widehat\sigma_t)\Delta(y).
\end{align*}
\eprop
\bew
i) For $x\in M$ we have $\Delta(x)=V(x\ot 1)V^*$ and then
\begin{align*}
\Delta(\sigma'_t(x))
&=V(\nabla^{it} x \nabla^{-it}\ot 1)V^* \\
&=V(\nabla^{it} x \nabla^{-it}\ot \widehat\nabla^{it}\widehat\nabla^{-it})V^* \\
&=(\nabla^{it}\ot \widehat\nabla^{it})V( x\ot 1)V^* (\nabla^{-it}\ot \widehat\nabla^{-it}) \\
&=(\sigma'_t\ot \tau_{-t})\Delta(x).
\end{align*}
\ssnl
ii) Similarly, for $y\in\widehat M$ we have $\widehat\Delta(y)=V^*(1\ot y)V$ and so
\begin{align*}
\Delta(\widehat\sigma_t(y))
&=V^*(1\ot \widehat\nabla^{it} y \widehat\nabla^{-it})V \\
&=V^*(\nabla^{it} \nabla^{-it}\ot \widehat\nabla^{it}y\widehat\nabla^{-it})V \\
&=(\nabla^{it}\ot \widehat\nabla^{it})V^*( 1\ot y)V (\nabla^{-it}\ot \widehat\nabla^{-it}) \\
&=(\widehat\tau_t\ot \widehat\sigma_t)\Delta(y).
\end{align*}
\ebew

Next we use the following formulas.

\prop
 \begin{equation*}
(\iota\ot\Delta)V=V_{12}V_{13}
\tussenen
(\widehat\Delta\ot 1)V^*=V^*_{13}V^*_{23}.
\end{equation*}
\eprop
\bew
i) For $a\in A$ and $\xi,\xi'\in\mathcal H$  we have
\begin{align*}
((\iota\ot\Delta)U)(\Lambda(a)\ot\xi\ot\xi')
&=\sum_{(a)}\Lambda(a_{(1)})\ot \Delta(a_{(2)}(\xi\ot\xi')\\
&=\sum_{(a)}\Lambda(a_{(1)})\ot a_{(2)}\xi\ot a_{(3)}\xi'\\
&=\sum_{(a)}U_{12}(\Lambda(a_{(1)})\ot \xi\ot a_{(2)}\xi')\\
&=U_{12}U_{13}(\Lambda(a)\ot \xi\ot\xi').
\end{align*}
ii) Similarly, using that $V^*(\xi\ot \widehat\Lambda(b))=\sum_{(b)} b_{(1)}\ot\widehat\Lambda(b_{(2)})$ we find the other equation.

\ebew
Then we obtain the following formulas.

\prop
For $x\in M$ and $y\in\widehat M$ we have
\begin{equation*}
\Delta(\tau_t(x))=(\tau_t\ot\tau_t)\Delta(x)
\tussenen
\widehat\Delta(\widehat\tau_t(y))=(\widehat\tau_t\ot\widehat\tau_t)\Delta(y)
\end{equation*}
\eprop
\bew
i) Multiply $(\iota\ot\Delta)V=V_{12}V_{13}$ with $\nabla^{it}\ot \widehat\nabla^{it}\ot \widehat\nabla^{it}$ on the left and with its inverse on the right. For the left and side we get
\begin{equation*}
(\iota\ot\tau_t\ot\tau_t)(\iota\ot \Delta)(V)
\end{equation*}
while on the right we get
\begin{align*}
(1\ot &\widehat\nabla^{it}\ot \widehat\nabla^{it}) 
(V_{12}V_{13})(1\ot \widehat\nabla^{-it}\ot \widehat\nabla^{-it})\\
&=(1\ot \widehat\nabla^{it}\ot \widehat\nabla^{it}) V_{12}(1\ot \widehat\nabla^{-it}\ot \widehat\nabla^{-it})\\
&\hspace*{40pt}(1\ot \widehat\nabla^{it}\ot \widehat\nabla^{it}) V_{13}(1\ot \widehat\nabla^{-it}\ot \widehat\nabla^{-it})\\
&=(\nabla^{-it}\ot 1\ot 1)V_{12}(\nabla^{it}\ot 1\ot 1)(\nabla^{-it}\ot 1\ot 1)V_{13}(\nabla^{it}\ot 1\ot 1)\\
&=(\nabla^{-it}\ot 1\ot 1)V_{12}V_{13}(\nabla^{it}\ot 1\ot 1)\\
&=(\nabla^{-it}\ot 1\ot 1)((\iota\ot\Delta)V)(\nabla^{it}\ot 1\ot 1)\\
&=(\iota\ot\Delta)(\nabla^{-it}\ot 1)V(\nabla^{it}\ot 1))\\
&=(\iota\ot\Delta)(1\ot\widehat\nabla^{it})V(1\ot \widehat\nabla^{-it})).
\end{align*}
Combining the two results we find that 
\begin{equation*}
(\iota\ot\tau_t\ot\tau_t)(\iota\ot \Delta)(V)
=(\iota\ot\Delta)((\iota\ot\tau_t)V)
\end{equation*}
and therefore $(\tau_t\ot\tau_t)\Delta(x)=\Delta(\tau_t(x))$ for all $x$.
\ssnl
In a completely similar way, we get from $(\widehat\Delta\ot 1)V^*=V^*_{13}V^*_{23}$ that $\widehat\Delta(\widehat\tau_t(y))=(\widehat\tau_t\ot\widehat\tau_t)\Delta(y)$ for all $y\in \widehat M$.
\ebew
\prop
For all $x\in M$ we have $\Delta(R(x))=\zeta(R\ot R)\Delta(x)$.
\eprop
\bew
By taking the adjoint of the equation $(\iota\ot\Delta(V)=V_{12}V_{13}$, it will easily follow from the first formula in the above proposition that $R$ flips the coproduct.
\ebew

We have seen that these maps leave $A$ invariant (see Proposition \ref{prop:2.13}). It follows that they also leave the C$^*$-algebra completion invariant.

\end{document}